\documentclass{amsproc}

\usepackage{amssymb, amsmath, amscd, latexsym, mathrsfs}
\usepackage{amsmath}
\usepackage{amsfonts}
\usepackage{appendix}
\usepackage[cmtip,all]{xy}

\newtheorem{thm}{Theorem}[subsection]
\newtheorem{lem}[thm]{Lemma}
\newtheorem{cor}[thm]{Corollary}
\newtheorem{prop}[thm]{Proposition}

\theoremstyle{definition}
\newtheorem{defn}[thm]{Definition}
\newtheorem{expl}[thm]{Example}

\theoremstyle{remark}
\newtheorem{rem}[thm]{Remark}
\numberwithin{equation}{subsection}  

\newcommand{\lra}{\longrightarrow}
\newcommand{\co}{\colon\!}

\newcommand{\smin}{\smallsetminus}

\newcommand{\id}{\textup{id}}

\newcommand{\holim}{\textup{holim}}
\newcommand{\hocolim}{\textup{hocolim}}

\newcommand{\hofiber}{\textup{hofiber}}
\newcommand{\mor}{\textup{mor}}
\newcommand{\map}{\textup{map}}
\newcommand{\fat}{\textup{fat}}

\newcommand{\emb}{\textup{emb}}

\newcommand{\man}{\mathsf{Man}}
\newcommand{\tw}{\textup{tw}}

\newcommand{\intr}{\textup{int}}

\newcommand{\sC}{\mathcal C}
\newcommand{\sD}{\mathcal D}

\newcommand{\sH}{\mathcal H}

\newcommand{\sM}{\mathcal M}

\newcommand{\sP}{\mathcal P}
\newcommand{\sQ}{\mathcal Q}

\newcommand{\sO}{\mathcal O}

\newcommand{\op}{\textup{op}}

\newcommand{\NN}{\mathbb N}
\newcommand{\RR}{\mathbb R}
\newcommand{\ZZ}{\mathbb Z}

\newcommand{\twosub}[2]{\begin{array}{cc}
\scriptstyle{#1} \\  [-1mm] \scriptstyle{#2}  \end{array}}

\newcommand{\threesub}[3]{\begin{array}{cc}
\scriptstyle{#1} \\  [-1mm] \scriptstyle{#2} \\ [-1mm] \scriptstyle{#3} \end{array}}

\newcommand{\colimsub}[1]{\begin{array}[t]{cc} \textup{colim} \\
[-1.7mm] \scriptstyle{#1} \end{array}}

\newcommand{\holimsub}[1]{\begin{array}[t]{cc} \textup{holim} \\ [-1mm]
\scriptstyle{#1} \end{array}}

\newcommand{\limsub}[1]{\begin{array}[t]{cc} \textup{lim} \\ [-1mm]
\scriptstyle{#1} \end{array}}

\newcommand{\hocolimsub}[1]{\begin{array}[t]{cc} \textup{hocolim} \\
[-1.7mm] \scriptstyle{#1} \end{array}}

\newcommand{\Tot}{\textup{Tot}}
\newcommand{\Tots}{\textup{Tots}}

\begin{document}

\title{Occupants in manifolds}
\author{Steffen Tillmann and Michael S. Weiss}%
\address{Math.~Institut, Universit\"at M\"unster, 48149 M\"unster, Einsteinstrasse 62, Germany}%
\email{s${}_-$till05@uni-muenster.de} \email{m.weiss@uni-muenster.de}





\begin{abstract} Let $K$ be a subset of a smooth manifold $M$. In some cases, functor calculus methods lead to
a homotopical formula for $M\smin K$ in terms of the spaces $M\smin S$, where $S$ runs through
the finite subsets of $K$.
\end{abstract}

\thanks{This project was supported by the Humboldt foundation through a Humboldt professorship for Michael Weiss, 2012-2017.}
\maketitle

\section{Occupants in a submanifold} \label{sec-sensub}
\subsection{Formulation of the problem}
Imagine a smooth manifold $M$ and a compact smooth submanifold $L$, both with empty boundary, of dimensions $m$ and $\ell$
respectively. We look for a homotopical description of $M\smin L$ in terms of the spaces $M\smin S$,
where $S$ runs through the finite subsets of $L$. The finite subsets $S$ of $L$ could be regarded as finite sets of
occupants.

\medskip
For one of the more geometric formulations of the problem, choose a Riemannian metric on $L$. Instead of working with finite
subsets $S$ of $L$, we work with thickenings of finite subsets of $L$ and we pay attention to inclusions of one such
thickening in another. More precisely
we work with pairs $(S,\rho)$ where $S$ is a finite subset of $L$ and $\rho$ is a function from $S$ to the positive
real numbers subject to two conditions.
\begin{itemize}
\item[-] For each $s\in S$, the exponential map $\exp_s$ at $s$ is defined and regular on the (compact) disk
of radius $\rho(s)$ about the origin in $T_sL$\,.
\item[-] The images in $L$ of these disks under the exponential
maps $\exp_s$ are pairwise disjoint.
\end{itemize}
For a pair $(S,\rho)$ satisfying the two conditions, let $V_L(S,\rho)\subset L$ be the union of the open balls of
radius $\rho(s)$ about points $s\in S$. Then $V_L(S,\rho)$ is an open subset of $L$ and is diffeomorphic to $\RR^\ell\times S$. The inclusion of
$M\smin V_L(S,\rho)$ in $M\smin S$ is a homotopy equivalence.

\smallskip
Let $C_k(L)$ be the space of unordered configurations of $k$ points in $L$.
For fixed $k\ge 0$, the pairs $(S,\rho)$ that satisfy the two conditions and the additional condition $|S|=k$ form
a space $C_k^\fat(L)$, an open subspace of the total space of some $k$-dimensional vector bundle on $C_k(L)$.
The forgetful projection $C_k^\fat(L)\to C_k(L)$ is a fiber bundle projection with contractible fibers, hence a homotopy equivalence.
Form the topological disjoint union of the spaces $C_k^\fat(L)$
for all $k\ge 0$ and view that as a (topological) poset $\sP(L)$ where
\[  (S,\rho) \le (T,\sigma) \]
means simply that $V_L(S,\rho)\subset V_L(T,\sigma)$. The poset $\sP(L)$ can also be viewed as a category. A contravariant
functor $\Phi$ from $\sP(L)$ to spaces is defined by
\begin{equation} \label{eqn-Phi}  \Phi(S,\rho)=M\smin V_L(S,\rho). \end{equation}
There is a map
\begin{equation} \label{eqn-sentry1}  M\smin L ~~\lra \holim~\Phi~, \end{equation}
determined by the inclusions $M\smin L\to \Phi(S,\rho)$ for $(S,\rho)\in \sP(L)$; the precise meaning
of $\holim~\Phi$ will be clarified in a moment.
With a view to occupants and the problem of finding unoccupied places, we ask:
\begin{itemize}
\item[-] is the map~(\ref{eqn-sentry1}) a weak homotopy equivalence~?
\end{itemize}
The question has a more numerical variant. For $j\ge 0$ let $\sP j(L)$ be the subspace and full topological sub-poset of $\sP(L)$
consisting of all $(S,\rho)$\ in $\sP(L)$ that satisfy $|S|\le j$. There is a map
\begin{equation} \label{eqn-sentry2}
M\smin L ~~\lra \holim~\Phi|_{\sP j(L)}
\end{equation}
determined by the inclusions $M\smin L\to \Phi(S,\rho)$. We ask:
\begin{itemize}
\item[-] is the map~(\ref{eqn-sentry2}) highly connected~? Is there a lower bound for the connectivity, expressed in terms of $j$~, which tends to
infinity with $j$~?
\end{itemize}
In reading these questions, keep in mind that
$\Phi$ has some continuity properties. This affects the meaning or definition of the
homotopy inverse limits. Here is a quick definition using the fact that
$\Phi$ is a subfunctor of a constant functor. More details can be found in definition~\ref{defn-bigholim} below.
Let $N\sP(L)$ be the nerve of $\sP(L)$,
a simplicial space. So $N_r\sP(L)$ is the \emph{space}
of order-reversing maps $u$ from $[r]=\{0,1,\dots,r\}$ to $\sP(L)$. \emph{Order-reversing}
means that $u(0)\ge u(1)\ge \cdots \ge u(r)$ in $\sP(L)$. Now
$\holim~\Phi$ in~(\ref{eqn-sentry1}) can be
described as a subspace of the space of all maps from the geometric realization $|N\sP(L)|$ to $M$, with the
compact-open topology.
A map $f\co |N\sP(L)|\to M$ belongs to that subspace if and only if for every $r$ and $u\in N_r\sP(L)$
with characteristic map $c_u\co \Delta^r\to |N\sP(L)|$, the composition $fc_u$ lands in $\Phi(u(r))\subset M$.
This description gives a rather good idea what $\holim~\Phi$ is: the space of all homotopy coherent ways to
choose a place in $M$ when some places in $L$ are occupied.
--- The homotopy limit in~(\ref{eqn-sentry2}) can be defined analogously.

\begin{thm}\label{thm-sentries} If $m-\ell\ge 3$, then the map
in~\emph{(\ref{eqn-sentry1})} is a weak homotopy equivalence and the map
in~\emph{(\ref{eqn-sentry2})} is $(1+(j+1)(m-\ell-2))$-connected.
\end{thm}

The proof of theorem~\ref{thm-sentries} is given at the end of this section. It is based on a reduction to standard
theorems in manifold calculus as found in \cite{WeissEmb} and \cite{GoWeEmb}. (See also remark~\ref{rem-oldvsnewcalc} below.)
The main idea is this: apply manifold calculus to the contravariant functor $F$ defined by $F(V)=M\smin V$
for open subsets $V$ of $L$. Then the left-hand side of~(\ref{eqn-sentry1}) is $F(L)$ and the right-hand side
is very reminiscent of $(T_\infty F)(L)$, where $T_\infty F$ is the ``Taylor series'' of $F$ in the sense of manifold calculus. Therefore the work
consists mainly in showing that $F$ is analytic. (But the definition $F(V):=M\smin V$ should be
regarded as provisional.)

\medskip
We begin with a slightly more systematic description or definition of the homotopy inverse
limits in~(\ref{eqn-sentry1}) and~(\ref{eqn-sentry2}). The idea is
to use the trusted formula of Bousfield-Kan \cite{BousfieldKan} while paying attention to topologies where appropriate.
Let $\Gamma_r(\Phi)$ be the space of (continuous) sections, with the compact-open topology,
of the fiber bundle $E_r^!(\Phi)\to N_r\sP(L)$ such that the fiber over a point $u\in N_r\sP(L)$ is $\Phi(u(r))$. So if
$u$ is given by a string $(S_0,\rho_0)\ge  \cdots\ge (S_r,\rho_r)$
in $\sP(L)$\,, then the fiber over $u$ is $\Phi(S_r,\rho_r)=M\smin V_L(S_r,\rho_r)$.

\begin{defn} \label{defn-bigholim}
{\rm The homotopy limit in~(\ref{eqn-sentry1}) is $\Tot\left([r]\mapsto \Gamma_r(\Phi)\right)$.
}
\end{defn}

\begin{rem} \label{rem-Reedy} The cosimplicial space $[r]\mapsto \Gamma_r(\Phi)$ is Reedy fibrant. This means that
for each $r\ge 0$ the \emph{matching map}
\[  \Gamma_r(\Phi)\lra \lim_{[r]\to[q]} \Gamma_q(\Phi) \]
determined by the non-identity (co)degeneracy operators is a fibration.
(Here $q<r$ and $[r]\to [q]$ stands for a monotone surjection;
see e.g.~\cite{Dugger} for more details.) It is a desirable property to have because
a map $X\to Y$ between Reedy fibrant cosimplicial spaces which is a degreewise weak
equivalence induces a weak equivalence $\Tot(X)\to \Tot(Y)$. Sketch of a proof showing that
$[r]\mapsto \Gamma_r(\Phi)$ is Reedy fibrant: since $E^!_r(\Phi)\to N_r\sP(L)$ is a fiber
bundle it is enough to note that the simplicial space $[r]\mapsto N_r\sP(L)$ is Reedy
cofibrant. (We use HELP, the homotopy extension lifting property.) It is Reedy cofibrant because, for every $r\ge 0$, the latching map
\[ \colimsub{[r]\to[q]} N_q\sP(L) \lra N_r\sP(L)  \]
(where $q<r$ etc.) is the inclusion of one ENR, Euclidean neighborhood retract, in another ENR as a
closed subspace. Such an inclusion is a cofibration \cite[ch.III,Thm.3.2]{Hu}. \newline
As a corollary of this observation, the homotopy type of $\holim~\Phi$ is independent of the choice
of Riemannian metric on $L$. If we choose two distinct Riemannian metrics $\mu$ and $\mu'$ on $L$
and call the corresponding functors $\Phi$ and $\Phi'$~, then there is a smooth path of Riemannian metrics
from $\mu$ to $\mu'$. Using that it is easy to produce a zigzag of degreewise weak equivalences
between Reedy fibrant cosimplicial spaces starting with the cosimplicial space $[r]\mapsto \Gamma_r(\Phi)$ and ending
with $[r]\mapsto \Gamma_r(\Phi')$. This leads to a zigzag of weak equivalences relating $\holim~\Phi$
to $\holim~\Phi'$. For descriptions of (something weakly equivalent to) $\holim~\Phi$ which do
not rely on any choice of Riemannian metric on $L$, see \cite{WeissConfHo} which is a sequel to this paper.
\end{rem}

\subsection{Discrete variants}
Write $\delta\sP(L)$ for the discrete variant of $\sP(L)$\,. So $\delta\sP(L)$ is a discrete poset and there
is a map of posets $\delta\sP(L)\to \sP(L)$ which is bijective and \emph{full}. That is, $(S,\rho)\le (T,\sigma)$ has the
same meaning in $\delta\sP(L)$ and in $\sP(L)$. That map from $\delta\sP(L)$ to $\sP(L)$ induces a map of homotopy inverse limits:
\begin{equation} \label{eqn-holimdiscrete}
\xymatrix@M=8pt{
{\rule{0mm}{6mm}\holimsub{(S,\rho)\textup{ in }\sP(L)} \Phi(S,\rho)} \ar[r] &
{\rule{0mm}{6mm}\holimsub{(S,\rho) \textup{ in }\delta\sP(L)} \Phi(S,\rho)}
}
\end{equation}
Similarly, restricting cardinalities of configurations we have a comparison map
\begin{equation} \label{eqn-holimdiscretek}
\xymatrix@M=8pt{
{\rule{0mm}{6mm}\holimsub{(S,\rho) \textup{ in }\sP j(L)} \Phi(S,\rho)} \ar[r] &
{\rule{0mm}{6mm}\holimsub{(S,\rho) \textup{ in }\delta\sP j(L)} \Phi(S,\rho)}
}
\end{equation}

\begin{lem} \label{lem-feed1} The maps~\emph{(\ref{eqn-holimdiscrete})} and~\emph{(\ref{eqn-holimdiscretek})} are weak
equivalences.
\end{lem}

\proof
We prove that the map~(\ref{eqn-holimdiscrete}) is a weak equivalence; the other statement
has a similar proof. ---
The first idea is to replace the topological poset $\sP(L)$ by a simplicial poset $[t]\mapsto P_t$\,.
Therefore let $P_t$ be the set (alias discrete space) of continuous maps from $\Delta^t$
to the underlying space of $\sP(L)$.
For $\sigma,\tau\in P_t$ we write $\sigma\le \tau$ to mean that
$\sigma(x)\le \tau(x)$ for all $x\in \Delta^t$. In particular the poset $P_0$ is identified with $\delta\sP(L)$.

For a morphism in $\Delta$ alias monotone map $\alpha\co[t]\to[u]$,  let $\Phi_\alpha$ be the contravariant functor
from $P_u\times\Delta^t$ to spaces obtained by composing
\[
\xymatrix@C=45pt{
P_u\times\Delta^t \ar[r]^{\id\times\alpha_*} &  P_u\times\Delta^u \ar[r]^-{\textup{evaluation}} & \sP(L)
}
\]
with $\Phi$.
Here $\Delta^t$ and $\Delta^u$ are regarded as topological posets with a trivial ordering but with a nontrivial (standard)
topology. In particular, $\Phi_\alpha$ is $\Phi|_{\delta\sP(L)}$ when $\alpha$ is the identity map of the object $[0]$ in $\Delta$.
Since $\Phi_\alpha$ is a continuous functor, we need to be explicit about $\holim~\Phi_\alpha$ if we want to use it.
It is defined as
\[ \Tot\Big([r]\mapsto \prod_{\sigma_0\ge\cdots\ge
\sigma_r~\in P_u} \lim(\Phi\circ\sigma_r\alpha_*)\Big) \]
where $\alpha_*\co \Delta^t\to\Delta^u$ is the map covariantly induced by
$\alpha$ and $\lim(\Phi\circ\sigma_r\alpha_*)$ is short for the space of maps $\Delta^t\to M$ which take $x\in\Delta^t$
to an element of
\[ \Phi(\sigma_r(\alpha_*(x)))=M\smin V_L(\sigma_r(\alpha_*(x)))\subset M\,. \]
Now $\alpha \mapsto \holim~\Phi_\alpha$
is a contravariant functor from the twisted arrow category $\tw(\Delta)$ of $\Delta$ to the category of spaces. (The objects
of the twisted arrow category $\tw(\sD)$ of a small category $\sD$ are morphisms $f\co c\to d$ in $\sD$,
and a morphism from $f\co c_0\to d_0$ to $g\co c_1\to d_1$
is a pair of morphisms $h\co c_0\to c_1$ and $k\co d_1\to d_0$ such that $f=kgh$.)
The map~(\ref{eqn-holimdiscrete}) can be obtained by composing the arrows in the lower row of the commutative diagram
\begin{equation}  \label{eqn-keydiscrete1}
\begin{split}
\xymatrix@M6pt{
&  {\limsub{\alpha\textup{ in }\tw(\Delta)} \holim~\Phi_\alpha} \ar[d]^-{\textup{incl.}} \ar[r]^-{\textup{incl.}} & {\prod_\alpha \holim~\Phi_\alpha}
\ar[d]^-{\textup{proj.}} \\
{\holim~\Phi} \ar[r] \ar[ur]  & {\holimsub{\alpha\textup{ in }\tw(\Delta)} \holim~\Phi_\alpha} \ar[r]^-{\textup{proj.}} & \holim~\Phi_{\id_0}
}
\end{split}
\end{equation}
where $\id_0$ is the identity morphism of $[0]$ in $\Delta$. (The map from $\holim~\Phi$
to $\prod_\alpha \holim~\Phi_\alpha$ in this diagram has coordinates equal to the prolongation maps
$\holim~\Phi\to \holim~\Phi_\alpha$.) We are going to show that the two arrows in the lower row
of diagram~(\ref{eqn-keydiscrete1}) are weak equivalences. For the one on the left it is a
routine task. We can write
\begin{align*} \holimsub{\alpha\co[t]\to[u]} \holim~\Phi_\alpha
&= &\holimsub{\alpha\co[t]\to[u]} \Tot\Big([r]\mapsto \prod_{\sigma_0\ge\cdots\ge
\sigma_r~\in P_u} \lim(\Phi\circ\sigma_r\alpha_*)\Big) \\
&=& \Tot\Big([r]\mapsto \holimsub{\alpha\co[t]\to[u]}\prod_{\sigma_0\ge\cdots\ge \sigma_r~\in P_u}\lim(\Phi\circ\sigma_r\alpha_*)\Big)
\end{align*}
where
\[  \holimsub{\alpha\co[t]\to[u]}\prod_{\sigma_0\ge\cdots\ge \sigma_r~\in P_u}\lim(\Phi\circ\sigma_r\alpha_*) \]
is the space of lifts as in the following commutative diagram:
\begin{equation} \label{eqn-keydiscrete2}
\begin{split}
\xymatrix@C=25pt@M=7pt{
 & E^!_r(\Phi) \ar[d] \\
{\rule{0mm}{6mm}\hocolimsub{\alpha\co[t]\to[u]} (N_rP_u\times\Delta^t)} \ar[r] \ar@{..>}[ur] \ar[dr]^-\simeq
& N_r\sP(L) \\
& \ar[u]^-{\textup{eval}}_-\simeq {\rule{0mm}{6mm}\colimsub{\alpha\co[t]\to[u]} (N_rP_u\times\Delta^t)}
}
\end{split}
\end{equation}
(Note that $N_rP_u$ is the set of singular $u$-simplices of the space $N_r\sP(L)$, and so the bottom term in diagram~(\ref{eqn-keydiscrete2})
is the geometric realization of the singular simplicial set of the space $N_r\sP(L)$.)
But since the horizontal map in diagram~(\ref{eqn-keydiscrete2})
is a weak equivalence, the map from the section space $\Gamma_r(\Phi)$ of $E_r^!(\Phi)\to N_r\sP(L)$
to
\[ \holimsub{\alpha\co[t]\to[u]}\prod_{\sigma_0\ge\cdots\ge \sigma_r~\in P_u}\lim(\Phi\circ\sigma_r\alpha_*) \]
that it induces is a weak equivalence for every $r\ge0$. It follows (by application of $\Tot$) that the first
arrow in the lower row of~(\ref{eqn-keydiscrete1}) is a weak equivalence. \newline
For the other arrow in the lower row of diagram~(\ref{eqn-keydiscrete1}) it suffices to show that the functor
\[  \big(\,\alpha\co[t]\to[u]\,\big)~\mapsto~\holim~\Phi_\alpha  \]
on $\tw(\Delta)$ takes all morphisms to weak equivalences.
This reduces easily to the weaker statement where we only consider morphisms between objects of the form $\alpha\co [0]\to[u]$
in $\tw(\Delta)$. For such an $\alpha$ it suffices to show that the map
\[ \holim~\Phi_{\id_0} \to \holim~\Phi_\alpha \]
induced by the unique morphism from $\alpha$ to $\id_0\co [0]\to[0]$ is a weak equivalence.
Now \cite[thm 6.12]{Dugger} can be applied. Then it only remains to show that the functor $\alpha^*\co P_u\to P_0$
is homotopically terminal, i.e., that for every element $z$ of $P_0$
the over category $P_{u,z}:=z/\alpha^*$ has a contractible classifying space. Note that
$P_{u,z}$ is just the full sub-poset of $P_u$ consisting of all $\sigma$ in $P_u$ such that
$z\le \sigma(w)$ in $P_0$\,, where $w$ is the vertex $\alpha_*(0)$ of $\Delta^u$. \newline
For $z=(S,\rho)\in P_0$ and $\sigma\in P_{u,z}$ let $F(\sigma)$ be the
space of \emph{based} continuous maps $g\co \Delta^u\to \emb(S,L)$ such that
$g(x)$ takes $S$ to $V_L(\sigma(x))$, for all $x\in \Delta^u$.
(\emph{Based} means that $g(w)$ is the inclusion $S\to L$.)
Clearly $F(\sigma)$ is contractible and $F$ is a covariant functor. So the projection
\[  \hocolim~F \lra |N(P_{u,z})| \]
is a weak equivalence. But there is also a forgetful projection
from $\hocolim~F$ to the space $Y$ of all based (continuous) maps $g\co \Delta^u\to \emb(S,L)$.
It is a fiber bundle projection with contractible base space.
Therefore it remains only to show that the fibers of this are contractible. Each fiber $Y_g$
has the form $|N(P_{u,z,g})|$ where $P_{u,z,g}$ is the full sub-poset of $P_{u,z}$ consisting of the $\sigma\in P_{u,z}$ which
in addition to the conditions for membership in $P_{u,z}$ satisfy $g(x)(S)\subset V_L(\sigma(x))$ for all $x\in \Delta^u$.
It is easy to see that the poset $P_{u,z,g}$ satisfies a form of directedness and so $|N(P_{u,z,g})|$ is
contractible. \qed

\begin{rem} \label{rem-holimdiscreteres} For an open subset $U$ of $L$, let $\sP(L)|_U$ be the full topological sub-poset of $\sP(L)$ consisting
of all $(S,\rho)$ such that $V_L(S,\rho)$ is contained in $U$. Let $\delta\sP(L)|_U$ be the corresponding discrete poset.
In the proof of theorem~\ref{thm-sentries}, we shall need a variant of lemma~\ref{lem-feed1} which states that the comparison map
\[
\xymatrix@M=8pt{
{\rule{0mm}{6mm}\holimsub{(S,\rho)\textup{ in }\sP(L)|_U} \Phi(S,\rho)} \ar[r] &
{\rule{0mm}{6mm}\holimsub{(S,\rho) \textup{ in }\delta\sP(L)|_U} \Phi(S,\rho)}
}
\]
is a weak equivalence. The proof is as for lemma~\ref{lem-feed1}.
\end{rem}

\subsection{Using manifold calculus}
Let $\sO(L)$ be the (discrete) poset of open subsets of $L$. In view of lemma~\ref{lem-feed1},
the following plan for a proof of theorem~\ref{thm-sentries} looks promising. There is a contravariant functor
\begin{equation} \label{eqn-funcompraw} V \mapsto M\smin V \end{equation}
from $\sO(L)$ to spaces. Manifold calculus as in \cite{WeissEmb} and \cite{GoWeEmb} was created to
help in understanding such functors. In particular, if a contravariant functor $F$ from $\sO(L)$ to
spaces has some reasonable properties such as isotopy invariance and satisfies some approximate excision conditions,
then manifold calculus has a formula
\[
\xymatrix{
F(L) \ar[r]^-{\simeq} & {\rule{0mm}{7mm}\holimsub{U\in\, \bigcup_k\!\!\sO k(L)} F(U)}~.
}
\]
(Here $\sO k(L)\subset \sO(L)$ is the full sub-poset whose elements are the open subsets of $L$ which are abstractly diffeomorphic to
$\RR^\ell\times S$ for some set $S$ with $|S|\le k$. Therefore
\[ \bigcup_k \sO k(L)\subset \sO(L) \]
is the full sub-poset whose
elements are the open subsets of $L$ which are abstractly diffeomorphic to $\RR^\ell\times S$ for some finite set $S$.)
Using lemma~\ref{lem-feed1}, we should be able to work from there to arrive at
theorem~\ref{thm-sentries}.
\newline
There is a small problem
with this plan. The functor~(\ref{eqn-funcompraw}) does not have all the good properties required
such as isotopy invariance. (Example: take $L$ to be $S^2=\RR^2\cup\infty$ and take $M$ to be $S^m=\RR^m\cup\infty$
for some $m\ge 2$. Let $V_1\subset\RR^2\subset L$
be the union of the open rectangles $]2^{-i-1},2^{-i}[\,\times\,]0,1[\,$ for $i=0,1,2,\dots$ and let $V_0\subset \RR^2\subset L$
be the union of the open squares $]2^{-i-1},2^{-i}[\,\times\,]0,2^{-i-1}[\,$ for $i=0,1,2,\dots$\,. The inclusion
$V_0\to V_1$ is an isotopy equivalence, but the induced homomorphism $\pi_{m-1}(M\smin V_1)\to \pi_{m-1}(M\smin V_0)$ is not
surjective.) But that is easy to fix. We rectify~(\ref{eqn-funcompraw}) by
setting
\begin{equation} \label{eqn-shapy} F(V): = \holimsub{C\subset V}~ M\smin C \end{equation}
for $V\in \sO(L)$, where $C$ runs through all compact subsets of $V$. Note that $F(L)$ has a forgetful projection map to $M\smin L$ which is
a weak equivalence, for if $V=L$ in~(\ref{eqn-shapy}), then there is a maximal choice for $C$ which is $C=L$. Moreover, that map $F(L)\to M\smin L$ has a
preferred section which can be obtained by composing $M\smin L\to \lim_C~(M\smin C)\to \holim_C~(M\smin C)$. The section is therefore also a weak equivalence.
Now we need to show that $F$ has reasonable properties such as isotopy invariance, and that it satisfies some
approximate excision conditions.

\begin{lem} \label{lem-good} The functor $F$ of~\emph{(\ref{eqn-shapy})} is \emph{good}. That is to say:
\begin{itemize}
\item[-] if $V_0\subset V_1$
are open subsets of $L$ such that the inclusion $V_0\to V_1$ is abstractly isotopic to a diffeomorphism, then the map
$F(V_1)\to F(V_0)$ induced by the inclusion is a weak homotopy equivalence;
\item[-] if $W\in \sO(L)$ is a union of open subsets $W_i$ where $i=0,1,2,3,\dots $ and $W_i\subset W_{i+1}$,
then the map from $F(W)$ to $\holim_i~F(W_i)$ determined by the inclusions $W_i\to W$ is a weak equivalence.
\end{itemize}
\end{lem}

\proof The second of the two properties claimed is obvious from the definition of $F$. By contrast the
first property is not easy to establish. Choose a sequence $(C_i)_{i\ge 0}$
of compact subsets of $V_1$ such that $C_{i}\subset C_{i+1}$ for all $i\ge 0$ and every compact subset of $V_1$
is contained in one of the $C_{i}$~. Then the projection from $F(V_1)$ to the sequential homotopy limit
\[  \holim_i~M\smin C_{i} \]
is a weak equivalence. Choose a smooth isotopy $(e_t\co V_0\to V_1)_{t\in[0,1]}$ such that
$e_0\co V_0\to V_1$ is the inclusion and $e_1\co V_0\to V_1$ is a diffeomorphism. (See remark~\ref{rem-isotopyinvariance}.)
Let $C_{t,i}:= e_t(e_1^{-1}(C_i))$ and note that $C_i$ has just been renamed $C_{1,i}$. Now the projection from $F(V_0)$ to
\[  \holim_i~M\smin C_{0,i} \]
is also a weak equivalence. Let $Y_i$ be the space of continuous maps $w\co [0,1]\to M$ such that
$w(t)\notin C_{t,i}$ for all $t\in [0,1]$. By a straightforward application of Thom's
isotopy extension theorem, the maps
\[ M\smin C_{0,i} \longleftarrow Y_i \lra M\smin C_{1,i}  \]
given by evaluation, $w\mapsto w(0)$ and $w\mapsto w(1)$, are homotopy equivalences. Therefore in the resulting
diagram of sequential homotopy limits
\[  \holim_i~M\smin C_{0,i} \longleftarrow \holim_i~ Y_i \lra \holim_i~M\smin C_{1,i}  \]
the two arrows are also weak equivalences. Now choose a monotone injective function
$\psi\co \NN \to \NN$ such that $C_{1,\psi(i)}$ contains $C_{t,i}$ for all $t\in [0,1]$.
This gives us a diagram
\[
\xymatrix@M=6pt@R=16pt@C=25pt{
F(V_1) \ar[rr]^-{\textup{induced by incl.}} \ar[d]^-\simeq && F(V_0) \ar[d]^-\simeq \\
\holim_i~M\smin C_{1,i} \ar[r]_-{\simeq}^-a & \holim_i~M\smin C_{1,\psi(i)} \ar[r]^b & \holim_i~M\smin C_{0,i} \\
&  \ar[ul]^-\simeq \holim_i~Y_i \ar[ur]_-\simeq
}
\]
where the vertical arrows and the arrow $a$ are forgetful projections while the arrow $b$
is induced by the inclusions $C_{0,i}\to C_{1,\psi(i)}$ for each $i\in\NN$.
The top square and the bottom triangle in the diagram are homotopy commutative; we leave the
verification to the reader. Therefore $b$ is a weak equivalence and then the top
horizontal arrow is a weak equivalence. \qed

\medskip
\proof[Proof of theorem~\ref{thm-sentries}]
We start by showing that the functor $F$ of~(\ref{eqn-shapy}) is analytic
and by giving some excision estimates for it. Since we know already that $F$ is good, it suffices to look into the following
situation. Let $P$ be a smooth compact codimension zero subobject (submanifold) of $L$ and let $Q_0,\dots Q_k$ be pairwise disjoint
compact codimension zero subobjects of $L\smin\intr(P)$. See \cite{GoWeEmb} for terminology. In more detail, this means that
$P$ is a compact smooth codimension zero submanifold with boundary of $L$ but each $Q_i$ is a compact smooth manifold with corners,
\[  \partial Q_i=\partial_0Q_i\cup \partial_1Q_i \]
where $\partial_0Q_i=Q_i\cap \partial P$ and $\partial_1Q_i$ is the closure of $\partial Q_i\smin P$ in $Q_i$\,.
There is the concept of handle index \cite[\S0]{GoWeEmb} of $Q_i$\,. (\emph{Important}: intuitively
this is the handle index of $Q_i$ relative to the subspace $\partial_0Q_i$\,.) We assume that
$Q_i$ has handle index $q_i$~. For $S\subset \{0,1,\dots,k\}$ let $Q_S:=P\cup\bigcup_{i\in S} Q_i$ and put
\[  W_S :=\intr(Q_S) \]
(taking the interior in $L$). The commutative cube of spaces
\[  S \mapsto F(W_S) \]
determines a map from $F(W_{\{0,1,\dots,k\}})$ to
\[  \holimsub{S\subsetneq \{0,1,2,\dots,k\}}~F(W_S)\,. \]
We need an estimate for the connectivity of that map, in terms of the dimensions $m$ and $\ell$,
the number $k$ and and the numbers $q_i$ for $i=0,1,\dots,k$. Because of our special assumption
we can work instead with the cube
\[  S ~~\mapsto~~G(W_S):= M\smin W_S \]
(contravariant in the variable $S$).
That cube is a strongly cocartesian cube in Goodwillie's terminology \cite{GoodwillieCalc2}. The map from
$G(W_{\{0,1,\dots,k\}})$ to $G(W_{\{0,1,\dots,k\}\smin\{i\}})$ is $(m-q_i-1)$-connected. Therefore the excision estimates of
\cite{EllisSteiner} and \cite[Thm.2.3]{GoodwillieCalc2} for such cubes apply here. Plugging these estimates
into \cite[Defn. 2.1]{GoWeEmb} we deduce that $F$ is $(m-2)$-analytic with excess $1$ for manifold calculus purposes. More to the point, the comparison map
\[   F(L) \lra \holimsub{U\in \,\,\bigcup_k\!\!\sO k(L)} F(U) \]
is a weak equivalence and the comparison map
\[ F(L) \lra \holimsub{U\in \,\sO k(L)} F(U) \]
is $(1+(k+1)(m-\ell-2))$-connected by \cite[Thm.~2.3]{GoWeEmb}. ---

Now we need to relate $\bigcup_k \sO k(L)$ to
$\delta\sP(L)$\,. There is a full monomorphism of posets
\[ \delta\sP(L)\to \bigcup_k \sO k(L) \]
which takes $(S,\rho)$ to $V_L(S,\rho)$.
This leads to a commutative diagram
\begin{equation} \label{eqn-round}
\begin{split}
\xymatrix@R=14pt@M=6pt{
M\smin L \ar[r] \ar[d]^-\simeq & {\rule{0mm}{6mm}\holimsub{(S,\rho)\in \delta\sP(L)} M\smin V_L(S,\rho)} \ar[d]_-\simeq \\
F(L) \ar[r]   & {\rule{0mm}{6mm}\holimsub{(S,\rho)\in \delta\sP(L)} F(V_L(S,\rho))} \\
F(L) \ar[u]_-= \ar[r]^-\simeq & {\rule{0mm}{6mm}\holimsub{U\in \,\,\bigcup_k\!\!\sO k(L)} F(U)} \ar[u]
}
\end{split}
\end{equation}
Since we want to know that the top horizontal arrow is a weak equivalence (keeping lemma~\ref{lem-feed1} in mind), we ought to show that
the lower right-hand vertical arrow, call it $g$, is a weak equivalence.
By a standard argument from the theory of homotopy limits \cite[thm.9.7]{DwyerKan} and with notation as in
remark~\ref{rem-holimdiscreteres}, the map $g$ is a weak equivalence if for each $U\in \,\,\bigcup_k\!\!\sO k(L)$
the map
\[  g_U\co F(U) \lra \holimsub{(S,\rho)\in \delta\sP(L)|_U} F(V_L(S,\rho)) \]
induced by the inclusions $V_L(S,\rho)\to U$ is a weak equivalence. To show that this is the case, fix $U$
and choose open $W\subset U$ such that $W=V_L(T,\sigma)$ for some $(T,\sigma)$ in $\sP(L)$ and the inclusion
$W\to U$ induces a bijection on $\pi_0$. Then in the commutative diagram
\[
\xymatrix{
F(U) \ar[d]^-{r_1} \ar[r]^-{g_U} & \ar[d]^-{r_2}  {\holimsub{(S,\rho)\in \delta\sP(L)|_U} F(V_L(S,\rho))} \\
F(W) \ar[r]^-{g_W} & {\holimsub{(S,\rho)\in \delta\sP(L)|_W} F(V_L(S,\rho))}
}
\]
the map $r_1$ is clearly a weak equivalence. The map $r_2$ is a weak equivalence, too.
Indeed by remark~\ref{rem-holimdiscreteres} it is allowed to replace $\delta\sP(L)|_U$ by $\sP(L)|_U$
and $\delta\sP(L)|_W$ by $\sP(L)|_W$, which means that in the right-hand column we can have homotopy
limits in the style of definition~\ref{defn-bigholim}. Then the verification amounts to seeing that a certain
map between cosimplicial spaces is a degreewise equivalence. Finally the
map $g_W$ is a weak equivalence by a cofinality argument. (There is a maximal element in
$\delta\sP(L)|_W$.) Therefore $g_U$ is a weak equivalence as claimed.

This takes care of the first part of theorem~\ref{thm-sentries}. The proof of the second part follows similar lines.
\qed

\begin{rem} The above proof of theorem~\ref{thm-sentries} might suggest that the $k$-th
Taylor approximation of the functor $F$ of~(\ref{eqn-shapy}), in the sense of manifold calculus,
can be obtained by post-composing $F$ with the $k$-th Taylor approximation of the identity functor from spaces to spaces,
in the sense of Goodwillie's homotopy functor calculus. Surprisingly, this is false. (It is also easy to
see that it is false in the case $k=1$.) A partial explanation is as follows. If $G$ is a functor from spaces
to spaces which is polynomial of degree $\le k$ in the sense of homotopy functor calculus,
then $GF$ is polynomial of degree $\le k$ in the manifold calculus sense. This is due to the similarity of the definitions of
\emph{polynomial functor} in the two functor calculuses, and to a property of $F$ which
was emphasized in the proof above. But if $G$ is a functor from spaces
to spaces which is homogeneous of degree $k$ in the sense of homotopy functor calculus,
then $GF$ need not be homogeneous of degree $k$ in the sense of manifold calculus. This is due to obvious differences in the
classification of homogeneous functors in the two functor calculuses.
\end{rem}

\begin{rem} \label{rem-failknot} How useful, interesting or faithful is the map~(\ref{eqn-sentry1})
when the codimension $m-\ell$ is less than 3~? Here is a codimension 2 case which is not encouraging. Let $M=S^3$
and let $L$ be a knot in $S^3$, your favorite knot, but not the unknot.
Let $\ZZ_\infty$ be the Bousfield-Kan $\ZZ$-completion functor from spaces to spaces. It comes with a
natural transformation $e\co \id\to \ZZ_\infty$. For simply connected spaces $X$, the natural
map $e\co X\to \ZZ_\infty X$ is a weak homotopy equivalence; this is applicable when $X$
is $\Phi(S,\rho)$ for some $(S,\rho)\in \sP(L)$. But for $X=M\smin L$ the natural map $X\to \ZZ_\infty X$ is not a weak
equivalence because $\ZZ_\infty(M\smin L)\simeq S^1$. (Instead it is a well-known map $M\smin L\to S^1$
which induces an isomorphism in ordinary integer homology.) We obtain a commutative diagram
\[
\xymatrix{
{M\smin L} \ar[d]^-e \ar[r]^-{(\ref{eqn-sentry1})} &  \ar[d]^-e_-\simeq {\holim~\Phi} \\
{\ZZ_\infty(M\smin L)} \ar[r] &  {\holim~(\ZZ_\infty\circ\Phi)}
}
\]
It follows that the map~(\ref{eqn-sentry1}), top horizontal arrow in the diagram, factors up to homotopy through the notorious
map $M\smin L\to S^1$. That seems to make~(\ref{eqn-sentry1}) tragically un-faithful, in this codimension 2 example.
\end{rem}

\begin{rem} \label{rem-isotopyinvariance} {\rm In articles on manifold calculus, the meaning of
\emph{isotopy equivalence} is sometimes ambiguous. According to one definition, call it (a), a smooth codimension zero embedding
$e\co U\to V$ (of smooth manifolds with empty boundary) is an isotopy equivalence if and only if there exists
an embedding $f\co V\to U$ such that $ef$ and $fe$ are smoothly isotopic to the respective identity maps.
According to another definition, call it (b), the embedding $e\co U\to V$ is an isotopy equivalence if and only if
it is isotopic (as a smooth embedding)
to a diffeomorphism from $U$ to $V$. We do not know whether definitions (a) and (b) are equivalent. Fortunately it is
easy to see that, if a functor from $\sO(L)$ to spaces takes isotopy equivalences as in definition (b) to weak equivalences,
then it takes isotopy equivalences as in definition (a) to weak equivalences. (The converse is obvious.) }
\end{rem}

\begin{rem} \label{rem-oldvsnewcalc} {\rm Two slightly different views exist on what manifold calculus is about.
In the older view laid out in \cite{WeissEmb} and \cite{GoWeEmb}, manifold calculus is about (some) contravariant functors
from $\sO(L)$ to spaces, where $L$ is a fixed background manifold. In a more modern view, described for example in
\cite{BoavidaWeiss} though it was also heralded in \cite{AroneTurchin}, manifold calculus is
about contravariant functors from a certain category $\man_\ell$ of \emph{all} smooth $\ell$-manifolds to spaces
(for some $\ell$). The morphisms in $\man_\ell$ are smooth embeddings between $\ell$-manifolds. More precisely, the morphisms
from $L_0$ to $L_1$ are organized into a space (or simplicial set), composition of morphisms is continuous (or is a simplicial map)
etc., which means that $\man_\ell$ is \emph{enriched} over the category of spaces (or simplicial sets). Similarly the category of
spaces is enriched over spaces (or simplicial sets), and the contravariant functors from $\man_\ell$ to spaces that we consider
in manifold calculus should respect the
enrichments.

A functor like $\emb(-,W)$ for a fixed smooth manifold $W$ lives comfortably in both settings:
the placeholder $-$ can be interpreted as an open subset of the fixed manifold $L$, or as an object of $\man_\ell$.
By contrast the functor $F$ of~(\ref{eqn-shapy}) which we have used in proving theorem~\ref{thm-sentries} seems to belong to the
older setting; $F(V)$ makes sense only for open subsets $V$ of $L$.

Does this mean that the modern reformulation of manifold calculus as in \cite{BoavidaWeiss} has thrown out the baby with the bathwater?
Nothing could be further from the truth. M.W. believes that most of the old manifold calculus can be subsumed in the new one as the branch
concerned with contravariant functors $G$ from $\man_\ell$ to spaces (preserving enrichment) which come equipped with a natural
transformation $\gamma$ to a representable functor
\[ \mor_{\man_\ell}(-,L)=\emb(-,L) \]
for fixed $L$ in $\man_\ell$\,. Such a pair $(G,\gamma)$ gives rise to a contravariant functor $G_\gamma$ from
$\sO(L)$ to spaces by
\[  G_\gamma(V) := \hofiber[\gamma\co G(V)\to \emb(V,L)] \]
for $V\in \sO(L)$, where the homotopy fiber is taken over the base point of $\emb(V,L)$. The construction $(G,\gamma)\to G_\gamma$
should be seen as a transform, i.e., it is often reversible. In particular
most of the contravariant functors from $\sO(L)$ to spaces that we encounter in the
old manifold calculus are weakly equivalent to $G_\gamma$ for some $G$ and $\gamma\co G\to\emb(-,L)$. \emph{Exercise}:
confirm this for the functor $F$ of~(\ref{eqn-shapy}). --- In any case, the proposed subsuming of the old manifold calculus in the
new one has not yet been carried out. That was one major reason for not using it here.
}
\end{rem}

\section{Occupants in the interior of a manifold} \label{sec-senabs}
\subsection{Formulation of the problem} \label{subsec-senabsformu}
Let $M$ be a smooth compact manifold with boundary. We look for a homotopical description
of $\partial M$ in terms of the spaces $M\smin S$, where $S$ runs through the finite subsets
of $M\smin \partial M$. To make that more precise, choose a Riemannian metric on $M$. Then $M\smin\partial M$ also
has a Riemannian metric and the topological poset $\sP(M\smin\partial M)$ is defined as in section~\ref{sec-sensub}.
Thus, elements of $\sP(M\smin\partial M)$ are pairs $(S,\rho)$
where $S$ is a finite subset of $M\smin\partial M$ and $\rho$ is a function from $S$ to the positive reals such that,
for each $s\in S$, the exponential map $\exp_s\co T_s(M\smin\partial M)\to M\smin\partial M$ is defined and regular on the disk of radius $\rho(s)$
about the origin, and the images of these disks in $M\smin\partial M$ are pairwise disjoint. For $(S,\rho)\in \sP(M\smin\partial M)$
let $V(S,\rho)\subset M\smin\partial M$ be the union of the open balls of radius $\rho(s)$ about points $s\in S$.
Then $V(S,\rho)$ is homeomorphic to $\RR^m\times S$.
A contravariant functor $\Psi$ from $\sP(M\smin\partial M)$ to spaces
is defined by
\begin{equation} \label{eqn-Psi}  \Psi(S,\rho)=M\smin V(S,\rho). \end{equation}
There are maps
\begin{equation} \label{eqn-sentry3}  \partial M ~\lra \holim~\Psi~, \end{equation}
\begin{equation} \label{eqn-sentry4}
\partial M~\lra \holim~\Psi|_{\sP j(M\smin\partial M)}
\end{equation}
induced by the inclusions $\partial M\to M\smin V(S,\rho)$.
The precise definition of the homotopy limits follows the pattern of definition~\ref{defn-bigholim}.

\begin{thm}\label{thm-sentriesabs} Suppose that $M$ is the total space of a smooth disk bundle $p\co M\to L$
of fiber dimension $c$ on a smooth closed manifold $L$. If $c\ge 3$, then the map~\emph{(\ref{eqn-sentry3})} is a weak equivalence
and the map~\emph{(\ref{eqn-sentry4})} is $(1+(j+1)(c-2))$-connected.
\end{thm}

The expression \emph{disk bundle} means just that: a smooth fiber bundle whose fibers are diffeomorphic to disks
$D^c$ of a fixed dimension $c$. It is not necessary to assume that $p$ is the disk bundle associated with a smooth vector
bundle on $L$.

\subsection{The tube lemma}
Theorem~\ref{thm-sentriesabs} will be proved by a reduction to theorem~\ref{thm-sentries}. The main idea for the reduction
is in lemma~\ref{lem-tube}
below. The lemma uses the notation of theorem~\ref{thm-sentriesabs}, but we can allow a disk bundle $p\co M\to L$ of any
fiber dimension $\ge 0$. We choose a Riemannian metric on $L$. Then in addition to the topological poset
$\sP(M\smin\partial M)$ there is the topological poset $\sP(L)$, and there are some interactions between the two which will be explored.
For $(S,\rho)\in \sP(M\smin\partial M)$ the open set $V(S,\rho)\subset M\smin\partial M$ was
defined just above. To be more consistent with
section~\ref{sec-sensub} we ought to write $V_{M\smin\partial M}(S,\rho)$,
but that would be cumbersome. For $(S,\rho)\in \sP(L)$ we still write $V_L(S,\rho)\subset L$ in the style
of section~\ref{sec-sensub}.

\begin{lem} \label{lem-tube}
The map
\begin{equation} \label{eqn-tube1} \hocolimsub{(S,\rho)\in\delta\sP(L)} C_k(p^{-1}(V_L(S,\rho))\smin\partial M)
~~\lra~~C_k(M\smin\partial M) \end{equation}
determined by the inclusions $C_k(p^{-1}(V_L(S,\rho))\smin\partial M)\to C_k(M\smin\partial M)$
is a weak equivalence.
\end{lem}

\proof It suffices to show that~(\ref{eqn-tube1}) is a Serre microfibration with
contractible fibers \cite[Lemma 2.2]{WeissHHA}. Contractibility of the fibers is straightforward. The fiber over a configuration
$T\in C_k(M\smin\partial M)$ is identified with the classifying space of the sub-poset $\sH_T$ of $\delta\sP(L)$
consisting of all $(S,\rho)$ such that $p(T)\subset V_L(S,\rho)$.
The poset $\sH_T$ is the target of a homotopy initial functor from the poset of the negative integers.

It remains to show that~(\ref{eqn-tube1}) is a Serre microfibration. Here
it is important to understand that although~(\ref{eqn-tube1}) and the projection to
$|N\delta\sP(L)|$ together determine an injective continuous map
\begin{equation} \label{eqn-tube2}
\hocolimsub{(S,\rho)\in\delta\sP(L)} C_k(p^{-1}(V_L(S,\rho))\smin\partial M)
\lra |N\delta\sP(L)|\times C_k(M\smin\partial M),
\end{equation}
that injective continuous map is not an embedding (homeomorphism onto the image).
Here is a lengthy example to illustrate the phenomenon and the advantages that it has for us.
Suppose for simplicity that $p\co M\to L$ is an identity map, i.e., disk bundle with fiber dimension $0$.
Take two elements $(S,\rho)$ and $(T,\sigma)$ of $\delta\sP(L)$ such that $(S,\rho)>(T,\sigma)$.
The inequality $(S,\rho)>(T,\sigma)$ determines a nondegenerate $1$-simplex in $N(\delta\sP(L))$, and an injective
map $\Delta^1\to |N\delta\sP(L)|$. View that as a path $w\co [0,1]\to N(\delta\sP(L))$, beginning at $(T,\sigma)$ and
ending at $(S,\rho)$. Suppose that $w$ has a lift $\bar w$ to a path in
$\hocolim_{(S,\rho)}~C_k(V_L(S,\rho))$. It is clear that the composition
\begin{equation} \label{eqn-digr}
\xymatrix{ [0,1] \ar[r]^-{\bar w} & {\rule{0mm}{6mm}\hocolimsub{(S,\rho)\in\delta\sP(L)} C_k(V_L(S,\rho))}  \ar[r] & C_k(L)
}
\end{equation}
has the form $t\mapsto R_t$ where the configuration $R_t$ is contained in
$V(T,\sigma)$ if $0\le t< 1$ and in $V(S,\rho)$ when $t=1$.
But more careful reasoning shows that
$R_1$ must be contained in $V(T,\sigma)$, too. (If a hint is needed, see remark~\ref{rem-taylor} below.)
This is fortunate for us because it implies at once that a
sufficiently small homotopy of the composition~(\ref{eqn-digr})
can be lifted to a homotopy of $\bar w$ itself, as the Serre microfibration condition wants to have it. ---
Now we return to our business, which is to establish the Serre microfibration condition
for the map~(\ref{eqn-tube1}).
Let $Z$ be a compact CW-space equipped with a map $f$ to the source in~(\ref{eqn-tube1}).
Using the injection~(\ref{eqn-tube2}), we may write $f=(f_1,f_2)$ where $f_1$ is a map with target
$|N\delta\sP(L)|$ and $f_2$ is a map with target $C_k(M\smin\partial M)$. Let
\[ h\co Z\times[0,1]\to C_k(M\smin\partial M) \]
be a homotopy such that $h_0=f_2$. Note that $f_1(Z)$
is contained in a finite union of cells of $|N\delta\sP(L)|$. For sufficiently small $\varepsilon>0$,
the map
\[  H\co Z\times[0,\varepsilon] \lra \rule{0mm}{6mm}\hocolimsub{(S,\rho)} C_k(p^{-1}(V_L(S,\rho))\smin\partial M) \]
defined by the formula $H=(f_1,h)$ using~(\ref{eqn-tube2}) is therefore continuous and well defined.
(The digression just above was meant to prepare for that slightly counter-intuitive claim.)
\qed

\begin{rem} \label{rem-taylor} On the theme of counter-intuitive topological properties of homotopy colimits in the
category of spaces, M.W. learned the following from Larry Taylor many years ago. Let $U$ be a bounded open interval in $\RR$. The mapping cylinder of
the inclusion $U\to\RR$ is not homeomorphic to a subspace of $\RR^2$. It is not metrizable.
\end{rem}

\begin{cor} \label{cor-tube} The map
\[ \hocolimsub{(S,\rho)\in\delta\sP(L)} N_0\sP(p^{-1}(V_L(S,\rho))\smin\partial M)
~~\lra~~N_0\sP(M\smin\partial M) \]
determined by the inclusions $p^{-1}(V_L(S,\rho))\smin\partial M\lra M\smin\partial M$
is a weak equivalence.
\end{cor}

\proof This is obtained from lemma~\ref{lem-tube} essentially by taking the disjoint union
over all $k\ge 0$, noting that the hocolim respects disjoint unions. Perhaps it should be
clarified that $\sP(U)$, for an open subset $U$ of $M\smin\partial M$, is defined or can be defined as the
full topological sub-poset of $\sP(M\smin\partial M)$ consisting of all elements $(S,\rho)$
such that $V(S,\rho)$ is contained in $U$ and has compact closure in $U$.
\qed

\begin{cor} \label{cor-tube2} For every $r\ge 0$ the map
\[ \hocolimsub{(S,\rho)\in\delta\sP(L)} N_r\sP(p^{-1}(V_L(S,\rho))\smin\partial M)
~~\lra~~N_r\sP(M\smin\partial M) \]
determined by the inclusions $p^{-1}(V_L(S,\rho))\smin\partial M\lra M\smin\partial M$
is a weak equivalence.
\end{cor}

\proof This is obtained from the previous corollary by noting that for open $U$ in $M\smin\partial M$
there is a homotopy pullback square
\[
\xymatrix{
N_r\sP(U) \ar[d] \ar[r] & N_r\sP(M\smin\partial M) \ar[d] \\
N_0\sP(U)  \ar[r] & N_0\sP(M\smin\partial M)
}
\]
where the horizontal arrows are inclusions and the vertical arrows are given by the ultimate target operator,
also known as 0-th vertex operator. (The square is also a strict pullback square, but this is less relevant for us.) \qed

\proof[Proof of theorem~\ref{thm-sentriesabs}, first part] There is a commutative square
\[
\xymatrix@M=6pt{
{\rule{0mm}{6mm}\holimsub{(T,\sigma)\in \sP(M\smin\partial M)}\Psi(T,\sigma)}
\ar[r] & {\rule{0mm}{9mm}\holimsub{(S,\rho)\in \delta\sP(L)} \!\!\!\!\!\!
\holimsub{\twosub{(T,\sigma)\in \sP(M\smin\partial M)}{\textup{cls. of }p\,(V(T,\sigma))\subset V_L(S,\rho)}} \Psi(T,\sigma)} \\
\partial M \ar[u] \ar[r] & {\rule{0mm}{6mm}\holimsub{(S,\rho)\in \delta\sP(L)} \partial M\cup\big(M\smin p^{-1}(V_L(S,\rho))\big)} \ar[u]
}
\]
By corollary~\ref{cor-tube2} the top horizontal arrow, given by specialization,
is a weak equivalence. By theorem~\ref{thm-sentries} and lemma~\ref{lem-feed1}, the lower horizontal arrow is a weak equivalence.
We want to know that the left-hand vertical arrow is a weak equivalence. So it suffices to show that
the right-hand vertical arrow is a weak equivalence. For that it suffices to show that
for fixed $(S,\rho)\in \delta\sP(L)$ the map
\[
\xymatrix@R=28pt{
{\holimsub{\twosub{(T,\sigma)\in \sP(M\smin\partial M)}{\textup{cls. of }p\,(V(T,\sigma))\subset V_L(S,\rho)}} \Psi(T,\sigma)} \\
{\rule{0mm}{6mm} \partial M\cup\big(M\smin p^{-1}(V_L(S,\rho))\big)} \ar[u]
}
\]
induced by the inclusion of $\partial M\cup\big(M\smin p^{-1}(V_L(S,\rho))\big)$ in the various
$\Psi(T,\sigma)$ is a weak equivalence. The target can also be written as
\begin{equation} \label{eqn-tub1} \holimsub{(T,\sigma)\in \sP(U)} \Psi(T,\sigma) \end{equation}
where $U$ is the open subset $p^{-1}(V_L(S,\rho))\smin\partial M$ of $M\smin\partial M$.
We now make a few alterations to that expression, which turn out to be weak equivalences under
$\partial M\cup\big(M\smin p^{-1}(V_L(S,\rho))\big)$.
\begin{itemize}
\item[(1)] Replace $\sP(U)$ by $\delta\sP(U)$.
\item[(2)] Let $F$ be the contravariant functor from $\sO(U)$ to spaces taking $W\in \sO(U)$ to
$\holim_{C}(M\smin C)$,
where $C$ runs through the compact subsets of $W$. Replace $\Psi(T,\sigma)$ by $F(V(T,\sigma))$.
\item[(3)] After implementing (1) and (2), replace $\delta\sP(U)$ by $\bigcup_k \sO k(U)$ and replace
$F(V(T,\sigma))$ for $(T,\sigma)\in \delta\sP(U)$ by $F(W)$ for $W\in \bigcup_k \sO k(U)$.
\end{itemize}
Alterations (1) and (3) can be justified by arguments which we have seen in section~\ref{sec-sensub}.
Alteration (2) is justified because there is a comparison map from $\Psi(T,\sigma)$
to $F(V(T,\sigma))$ which is a weak equivalence.
In this way, expression~(\ref{eqn-tub1}) turns into
\begin{equation} \label{eqn-tub2} \holimsub{W\in\,\,\bigcup_k\!\sO k(U)} F(W)\,. \end{equation}
But the poset $\bigcup_k \sO k(U)$ has a maximal element, which is $U$ itself. Therefore
expression~(\ref{eqn-tub2}) can be replaced by $F(U)$. It is easy to see that the
reference map from $\partial M\cup\big(M\smin p^{-1}(V_L(S,\rho))\big)$ to $F(U)$ is a weak
equivalence. \qed

\medskip
\proof[Proof of theorem~\ref{thm-sentriesabs}, second part] Fix $j>0$ as in~(\ref{eqn-sentry4}).
We need a modification of lemma~\ref{lem-tube}. Let $k$ be another integer such that $j\ge k\ge 0$.
The modification states that the projection map
\[\hocolimsub{(S,\rho)\in\delta\sP j(L)} C_k(p^{-1}(V_L(S,\rho))\smin\partial M)
~~\lra~~C_k(M\smin\partial M) \]
is a weak equivalence. The proof is exactly like the proof of lemma~\ref{lem-tube} itself: the map is again a
Serre microfibration with contractible fibers. We need $j\ge k$ for the contractibility of the fibers. ---
There is a commutative square
\[
\xymatrix@M=6pt@R=17pt{
{\rule{0mm}{6mm}\holimsub{(T,\sigma)\in \sP j(M\smin\partial M)}\Psi(T,\sigma)}
\ar[r] & {\rule{0mm}{9mm}\holimsub{(S,\rho)\in \delta\sP j(L)} \!\!\!\!\!\!
\holimsub{\twosub{(T,\sigma)\in \sP j(M\smin\partial M)}{\textup{cls. of }p\,(V(T,\sigma))\subset V_L(S,\rho)}} \Psi(T,\sigma)} \\
\partial M \ar[u] \ar[r] & {\rule{0mm}{6mm}\holimsub{(S,\rho)\in \delta\sP j(L)} \partial M\cup\big(M\smin p^{-1}(V_L(S,\rho))\big)} \ar[u]
}
\]
By a modification of corollary~\ref{cor-tube2} which is a consequence of the modification of lemma~\ref{lem-tube}
just formulated, the top horizontal arrow
is a weak equivalence. By theorem~\ref{thm-sentries} and lemma~\ref{lem-feed1}, the lower horizontal arrow is
$(1+(j+1)(c-2))$-connected.
We want to know that the left-hand vertical arrow is $(1+(j+1)(c-2))$-connected. So it suffices to show that
the right-hand vertical arrow is a weak equivalence. This can be verified as in the proof of the first half
of theorem~\ref{thm-sentriesabs}. \qed

\section{Gates} \label{sec-gate}
This section generalizes the previous two. Consequently it has two slightly different themes.

\subsection{Submanifold case}
For the first theme, imagine a smooth manifold $M$ with boundary and a neat smooth compact submanifold $L$, so that
$\partial L\subset \partial M$. We look for a homotopical description of $M\smin L$ in terms of the
spaces $M\smin S$, where $S$ runs through the finite subsets of $L\smin \partial L$.
In the case where $\partial L$ and $\partial M$ are empty, this is exactly the
situation of section~\ref{sec-sensub}. Also, in the case where $\partial L$ is empty but $\partial M$
is nonempty, it is almost exactly the situation of section~\ref{sec-sensub} because in
such a case it makes no substantial difference if we delete $\partial M$ from $M$. \newline
For a more precise formulation we extend the definition of $\sP(L)$ given in section~\ref{sec-sensub} so
that $L$ is allowed to have a nonempty boundary. Choose a Riemannian metric on $L$. The elements
of $\sP(L)$ are going to be pairs $(S,\rho)$ where $S$ is a finite subset of $L\smin\partial L$ and $\rho$ is a
function from $S\sqcup \partial L$ to the positive reals, locally constant
on $\partial L$ and subject to a few more conditions.
\begin{itemize}
\item[-] For each $s\in S$, the exponential map $\exp_s$ at $s$ is defined and regular on the disk
of radius $\rho(s)$ about the origin in $T_sL$\,.
\item[-] The (boundary-normal) exponential map is defined and regular on the set of all tangent vectors $v\in T_zL$ where $z\in \partial L$,
where the vector $v$ is inward perpendicular to $T_z\partial L$ and $|v|\le \rho(z)$.
\item[-] The images in $L$ of these disks and the image
of this band under the exponential map(s) are pairwise disjoint.
\end{itemize}
For a pair $(S,\rho)$ satisfying these conditions, let $V_L(S,\rho)\subset L$ be the union of the
open balls of radius $\rho(s)$
about elements $s\in S$ and the open collar on $\partial L$ determined the normal distance function $\rho_{|\partial L}$.
Then $V_L(S,\rho)$ is diffeomorphic
to $(\RR^\ell\times S)\sqcup [0,1[\,\times\partial L$ and the inclusion of $M\smin V_L(S,\rho)$ in
$M\smin S$ is a homotopy equivalence. The partial order on $\sP(L)$ is defined so that
$(S_0,\rho_0)\le (S_1,\rho_1)$ if and only if $V_L(S_0,\rho_0)\subset V_L(S_1,\rho_1)$.
In this partial order the boundary $\partial L$ acts like a gate which allows occupants to leave. \newline
The poset $\sP(L)$ can also be viewed as a category. A contravariant
functor $\Phi$ from $\sP(L)$ to spaces is defined by
\begin{equation} \label{eqn-Phibd}  \Phi(S,\rho)=M\smin V_L(S,\rho). \end{equation}
There is a map
\begin{equation} \label{eqn-sentry5}  M\smin L ~~\lra \holim~\Phi~, \end{equation}
determined by the inclusions $M\smin L\to \Phi(S,\rho)$ for $(S,\rho)\in \sP(L)$.
Also, let $\sP j(L)$ be the subspace and full topological sub-poset of $\sP(L)$
consisting of all $(S,\rho)$\ in $\sP(L)$ that satisfy $|S|\le j$. Then again there is a map
\begin{equation} \label{eqn-sentry6}
M\smin L ~~\lra \holim~\Phi|_{\sP j(L)}
\end{equation}
determined by the inclusions $M\smin L\to \Phi(S,\rho)$.

\begin{thm}\label{thm-sentrieswall} If $m-\ell\ge 3$, then the map~\emph{(\ref{eqn-sentry5})} is a weak
homotopy equivalence and the map~\emph{(\ref{eqn-sentry6})} is $(1+(j+1)(m-\ell-2))$-connected.
\end{thm}

The proof of this is very similar to the proof of theorem~\ref{thm-sentries} and the
details are therefore omitted.

\subsection{Absolute case} \label{subsec-absgate}
Our second topic is a generalization of theorem~\ref{thm-sentriesabs} to a situation with more complicated boundary
conditions. Let $M$ be a compact smooth manifold with boundary and corners in the boundary. In particular
$\partial M$ is the union of two codimension zero smooth submanifolds $\partial_0M$ and $\partial_1M$ that intersect in
the corner set
\[  \partial_0M\cap\partial_1M=\partial(\partial_0M)=\partial(\partial_1M). \]
We look for a homotopical description of $\partial_1M$ in terms of the spaces $M\smin S$, where $S$
runs through the finite subsets of $M\smin\partial M$. \newline
Now $M\smin\partial_1M$ is a smooth manifold with boundary $\partial_0M\smin\partial_1M$; both
$M\smin\partial_1 M$ and its boundary can be noncompact. Again we choose a Riemannian metric on all of $M$. For simplicity
we require it to be a product metric in a neighborhood of $\partial_1M$, i.e., the product of a Riemannian metric
on $\partial_1M$
and a Riemannian metric on $[0,\varepsilon]$ for some $\epsilon>0$.  
Then we can define a topological poset $\sP(M\smin\partial_1 M)$ roughly as in section~\ref{sec-senabs}.
The elements are pairs $(S,\rho)$ where $S$ is a finite subset of $M\smin\partial M=(M\smin\partial_1M)\smin\partial(M\smin\partial_1M)$
and $\rho$ is a function from $S\sqcup(\partial_0M\smin\partial_1M)$ to the positive reals which is locally constant on
$\partial_0M\smin\partial_1M$ and
subject to the usual conditions. The usual conditions are, briefly stated: regularity of the exponential maps on the tangential disks
and the tangential closed band defined by $\rho$, and pairwise disjointness of their images in $M\smin\partial_1 M$. The union
of the open balls of radius $\rho(s)$ about elements $s\in S$ and of the (half-)open band determined by the normal distance function
$\rho$ on $\partial_0M$ is denoted by $V(S,\rho)$. Therefore $V(S,\rho)$ is an open subset of $M\smin\partial_1M$ diffeomorphic
to $(\RR^m\times S)\sqcup(\partial_0M\smin\partial_1M)\times[0,1[\,$.

\smallskip
A contravariant functor $\Psi$ from $\sP(M\smin\partial_1 M)$ to spaces
is defined by
\begin{equation} \label{eqn-Psibd}  \Psi(S,\rho)=M\smin V(S,\rho). \end{equation}
There are maps
\begin{equation} \label{eqn-sentry7}  \partial_1 M ~\lra \holim~\Psi~, \end{equation}
\begin{equation} \label{eqn-sentry8}
\partial_1 M~\lra \holim~\Psi|_{\sP j(M\smin\partial M)}
\end{equation}
induced by the inclusions $\partial_1 M\to M\smin V(S,\rho)$.

\begin{thm} \label{thm-sentriesabsbd} Suppose that $M\smin\partial_1M$ has a neat compact smooth submanifold $L$ of codimension $c\ge 3$
making $M$ into a smooth thickening of $L\cup\partial_0M$~; see definition~\emph{\ref{defn-thick}} below.
Then the map~\emph{(\ref{eqn-sentry7})}
is a weak equivalence and the map~\emph{(\ref{eqn-sentry8})} is $(1+(j+1)(c-2))$-connected.
\end{thm}

\begin{defn} \label{defn-thick} {\rm On the meaning of \emph{smooth thickening} in theorem~\ref{thm-sentriesabsbd}:
the main points are that the inclusion of $\partial_0M\cup L$ in $M$ is a homotopy equivalence and
that the inclusion of $\partial_1M$ in $M\smin L$ is a homotopy equivalence. \newline
In detail, we start with  a smooth compact manifold $L$, another smooth compact
manifold $A$ (think $A=\partial_0M$) and a smooth embedding $u\co \partial L\to A$
which avoids the boundary of $A$. Then we can speak of $L\cup A$, the
pushout of $L\leftarrow \partial L\to A$. Let $Q\subset L$ be a closed collar, so that $\partial L\subset Q$
and there is a diffeomorphism $Q\to \partial L\times[0,1]$ extending the map $x\mapsto (x,0)$
on $\partial L$. Let $L_1$ be the closure of $L\smin Q$ in $L$.
To make a smooth thickening of $L\cup A$ we need in addition
\begin{itemize}
\item a smooth disk bundle
$E\to L_1$ whose total space has dimension $\dim(A)+1$;
\item an identification of $E|_{\partial L_1}$ with the normal disk bundle
of the embedding $u\co\partial L\to A$, over the evident diffeomorphism $\partial L_1\cong \partial L$.
\end{itemize}
Then the pushout $T$ of $A\times[0,1]\leftarrow E|_{\partial L_1} \rightarrow E$ is defined, on the understanding that the
left-hand arrow embeds $E|_{\partial L_1}$ into $A\times\{1\}$. Now $T$ is a compact manifold with boundary, smooth
with corners. The corner set has three disjoint parts: $\partial A\times\{0\}$, $\partial A\times\{1\}$ and
$\partial(E|_{\partial L_1})$.
Think of $\partial T$ as the union of $\partial_0T:=A\times\{0\}$ and $\partial_1T$, the closure
of $\partial T\smin \partial_0T$ in $\partial T$. The parts of the corner set not accounted for by $\partial_0T\cap \partial_1T$
should be subjected to smoothing. There is a copy of $L=Q\cup L_1$ contained in $T$. And of course there is also a copy of
$A$ contained in $T$, in the shape of $A\times\{0\}$.
Any smooth manifold with corners which is diffeomorphic to such a $T$ relative to $L\cup A=L\cup\partial_0T$ can be called
a smooth thickening of $L\cup A$ or of $L\cup\partial_0T$.
}
\end{defn}

\begin{expl} \label{rem-handsfree} In theorem~\ref{thm-sentriesabsbd}, the manifold $M$ can be
the total space of a smooth disk bundle $M\to L$ of fiber dimension $c\ge 3$, where
$L$ is allowed to be compact smooth \emph{with boundary}, and $M_0:= M|_{\partial L}$.
The proof of theorem~\ref{thm-sentriesabsbd} in this case is a straightforward variation of the proof of
theorem~\ref{thm-sentriesabs}.
\end{expl}

\subsection{Something like engulfing} \label{subsec-gategulf}
We turn to the proof of theorem~\ref{thm-sentriesabsbd}. This is broken up into remarks,
definitions, lemmas and even a corollary.

\begin{rem} \label{rem-onlysecondpart} The second part of theorem~\ref{thm-sentriesabsbd} (the high connectivity statement)
implies the first part (the weak equivalence statement). This is easy to see if we use the
definition of $\holim~\Psi$ as a subspace of the space of maps from $|N\sP(M\smin\partial_1M)|$
to $M$. The inclusion of $|N\sP j(M\smin\partial_1M)|$
in $|N\sP k(M\smin\partial_1M)|$, for $k>j$, is a cofibration. It follows that
the projection
\[  \holim~\Psi|_{\sP k(M\smin\partial_1M)} \lra  \holim~\Psi|_{\sP j(M\smin\partial_1M)} \]
is a fibration for $k>j$. Therefore the canonical inclusion
\[
\xymatrix@M=8pt{
 \holim~\Psi~=~\lim_j~\holim~\Psi|_{\sP j(M\smin\partial_1M)}
\ar[r] & \holim_j~\holim~\Psi|_{\sP j(M\smin\partial_1M)}
}
\]
is a weak equivalence. The homotopy groups $\pi_r$ of the right-hand side can be calculated
as inverse limits
\[ \lim_j \pi_r\left(\holim~\Psi|_{\sP j(M\smin\partial_1M)}\right). \]
The higher derived inverse limit $\lim^1$ does not contribute to this calculation
because of the Mittag-Leffler criterion \cite[ch.7,App.]{Switzer}. The criterion is applicable
here because we are assuming the second part
of theorem~\ref{thm-sentriesabsbd}.
\end{rem}

\begin{rem} Since the validity of theorem~\ref{thm-sentriesabsbd} does not depend on the Riemannian
metric which we select for $M$, we can choose a Riemannian metric with very convenient properties.
We shall assume that it is a product metric near $\partial_0M$ as well, i.e., a neighborhood
of $\partial_0M$ is isomorphic as a Riemannian manifold to a product $\partial_0M\times[0,\varepsilon]$
where the interval $[0,\varepsilon]$ has the standard Riemannian metric, and the isomorphism
takes $\partial_0M\subset M$ to $\partial_0M\times\{0\}\subset\partial_0M\times[0,\varepsilon]$.
(The product structure near $\partial_0M$
is automatically compatible with the product structure near $\partial_1M$ which we assumed earlier,
so that $\partial\partial_0M=\partial_0M\cap\partial_1M$ has a neighborhood in $M$ which is isomorphic to
$\partial\partial_0M\times[0,\varepsilon]\times[0,\varepsilon]$ as a Riemannian manifold.)
As a result we have a standard
compact collar for $\partial_0M$ (of width $\varepsilon$). We shall also assume that the intersection of
$L$ with that collar has the form $\partial L\times [0,\varepsilon]$ in the collar coordinates. \newline
In addition we choose an open tubular neighborhood $U$ of $L$ such that the closure $\bar U$ of $U$ in $M$
is a smooth disk bundle over $L$. Also, the intersection of $U$ with the standard collar on $\partial_0M$
(of width $\varepsilon$) is required to have the form $\partial U\times[0,\varepsilon]$ in the collar
coordinates. Here $\partial U\subset \partial_0M$ is a tubular neighborhood of $\partial L$ whose
closure in $\partial_0M$ is a smooth disk bundle over $\partial L$.
\end{rem}

\begin{defn} Let $\sC_1$ be the full topological sub-poset of $\sP j(M\smin \partial_1 M)$ consisting of the objects
$(S,\rho)$ such that the locally constant function $\rho|_{\partial_0M}$ is $\le \varepsilon$ everywhere
and the disks of radius $\rho(s)$ about elements $s\in S$ are all contained in $U$, the specified
tubular neighborhood of $L$. Write $\delta\sC_1$ for the discrete variant.
\end{defn}

\begin{defn}
Let $\sC_0$ be the full topological sub-poset of $\sP j(M\smin \partial_1 M)$ consisting of the objects $(S,\rho)$
such that $\rho\le \varepsilon/3j$, and the closure of $V(S,\rho)$ is contained in the union of $U$ and the
standard collar on $\partial_0M$ of width $\varepsilon/3j$.  Nota bene: for an object $(S,\rho)$ of $\sC_0$
it can happen that $S$ is not contained in $U$. Write $\delta\sC_0$ for the discrete variant.
\end{defn}

\begin{lem} \label{lem-funnycofinal}
For every element $(S,\rho)$ of $\delta\sC_0$ there is some element $(T,\sigma)$ of $\delta\sC_1$
such that $(T,\sigma)\ge (S,\rho)$ in $\sP j(M\smin\partial_1M)$. Indeed the
sub-poset of $\delta\sC_1$ consisting of the elements $(T,\sigma)$ which satisfy $(T,\sigma)\ge(S,\rho)$
has a contractible classifying space.
\end{lem}

\proof It is enough to note that there is a real number $\tau$, strictly between $\varepsilon/3j$ and
$\varepsilon$, such that the parallel hypersurface to $\partial_0M$ in $M$ at distance $\tau$ from $\partial_0M$
has empty intersection with the closure of $V(S,\rho)$. This is due to our assumption $|S|\le j$
and the smallness of the radii in the metric balls which are part of $V(S,\rho)$.
\qed

\begin{cor} \label{cor-funnycofinal} There is a homotopy commutative diagram of the following shape,
\[
\xymatrix{
 \holim~\Psi|_{\sP j(M\smin\partial_1M)} \ar[d] \ar[r] & \ar@{..>}[d] \holim~\Psi|_{\delta\sC_1} \\
 \holim~\Psi|_{\delta\sC_0} \ar@{..>}[r]^-{\simeq} & Y
}
\]
where the solid arrows are given by restriction.
\end{cor}

\proof Define $Y$ to be
\[
\holimsub{(S,\rho)\in \delta\sC_0} \holimsub{\threesub{(S,\rho)\ge (S'\!,\,\rho')\le (T,\sigma)}{(T,\sigma)\in \delta\sC_1}
{(S'\!,\,\rho')\in \delta\sC_0}} \Psi(S',\rho')\,.
\]
The inner homotopy limit is taken over the subset $\sM(S,\rho)$ of the product poset $\delta\sC_0\times\delta\sC_1$
consisting of pairs $((S'\!,\,\rho'),(T,\sigma))$ which satisfy $(S,\rho)\ge (S'\!,\,\rho')$
and $(S'\!,\,\rho')\le (T,\sigma)$ in $\sP j(M\smin \partial_1M)$. We note that $\sM(S,\rho)$ is covariantly
functorial in $(S,\rho)\in \delta\sC_0$~, so that $(S,\rho)\le (\bar S,\bar\rho)$ in $\delta\sC_0$ implies
$\sM(S,\rho)\subset \sM(\bar S,\bar\rho)$. This implies that the inner homotopy limit is a contravariant
functor of the variable $(S,\rho)$. --- The two broken arrows are then fairly obvious and the homotopy
commutativity of the resulting diagram is also clear. The horizontal
broken arrow is a weak equivalence by lemma~\ref{lem-funnycofinal}. More precisely, there are comparison maps
\[  \Psi(S,\rho) \quad\lra \holimsub{\twosub{(S,\rho)\ge (S'\!,\,\rho')}
{(S'\!,\,\rho')\in \delta\sC_0}} \Psi(S',\rho') \quad\lra \holimsub{\threesub{(S,\rho)\ge (S'\!,\,\rho')\le (T,\sigma)}{(T,\sigma)\in \delta\sC_1}
{(S'\!,\,\rho')\in \delta\sC_0}} \Psi(S',\rho').\]
the second of which is prolongation along a forgetful functor. The first of these is a weak equivalence because the indexing poset
in the target expression has a terminal (maximal) element. The second one is a weak equivalence because the forgetful functor
from $\sM(S,\rho)$ to the poset $\{\,(S',\rho')\in\delta\sC_0~|~(S,\rho)\ge(S',\rho')\}$ is homotopy terminal \cite{Dugger}
by lemma~\ref{lem-funnycofinal}.
Therefore the composition is a weak equivalence, and this is what we need to conclude that the lower horizontal arrow in the diagram
is a weak equivalence. \qed

\begin{lem} \label{lem-res1} The map $\holim~\Psi|_{\sP j(M\smin\partial_1M)}\lra \holim~\Psi|_{\sC_0}$
given by restriction is a weak equivalence. The map $\holim~\Psi|_{\sC_0}\lra \holim~\Psi|_{\delta\sC_0}$ given by
restriction is also a weak equivalence
\end{lem}

\proof For the first statement, observe that the inclusion map from $N_r\sC_0$ to $N_r\sP j(M\smin\partial_1M)$ is a
homotopy equivalence for every $r\ge 0$. The second statement can be proved like lemma~\ref{lem-feed1}. \qed

\begin{lem} \label{lem-hc} The composition
\[
\xymatrix@M=8pt{
{\partial_1M} \ar[r]^-{(\ref{eqn-sentry8})} & {\holim~\Psi|_{\sP j(M\smin\partial_1M)}} \ar[r] & {\holim~\Psi|_{\sC_1}}
}
\]
is $(1+(j+1)(c-2))$-connected.
\end{lem}

\proof Following the lines of section~\ref{sec-senabs}, reduce to the claim that the map
\begin{equation} \label{eqn-sigh}  \partial_1M \lra \holimsub{(S,\rho)\in \sP j(L)} M\smin V_L(S,\rho) \end{equation}
induced by the inclusions $\partial_1M\to M\smin V_L(S,\rho)$ is so-and-so highly connected. (Following the lines of section~\ref{sec-senabs}
gives a zigzag of weak equivalences \emph{under} $\partial_1M$
relating $\holim~\Psi|_{\sC_1}$ to the target in~(\ref{eqn-sigh}).) The map~(\ref{eqn-sigh})
can be written as a composition
\[ \partial_1M \hookrightarrow M\smin L \lra  \holimsub{(S,\rho)\in \sP j(L)} M\smin V_L(S,\rho). \]
The first arrow is a homotopy equivalence by our definition of \emph{thickening} and the
second arrow is so-and-so-highly connected by theorem~\ref{thm-sentrieswall}. \qed

\medskip
\proof[Conclusion of proof of theorem~\ref{thm-sentriesabsbd}] By remark~\ref{rem-onlysecondpart} we can
concentrate on the high connectivity statement, for a fixed $j\ge 0$.
The homotopy commutative square of corollary~\ref{cor-funnycofinal}
can be enlarged to a diagram
\[
\xymatrix@C=35pt@R=20pt{ \partial_1M \ar[r]^-{(\ref{eqn-sentry8})}  &
 \holim~\Psi|_{\sP j(M\smin\partial_1M)} \ar[dd]^-\simeq \ar[r] & \ar[d] \holim~\Psi|_{\sC_1} \\
    && \holim~\Psi|_{\delta\sC_1} \ar[d] \\
 & \holim~\Psi|_{\delta\sC_0} \ar[r]^-\simeq &  Y
}
\]
The weak equivalence labels in the diagram
are justified by lemma~\ref{lem-res1} and corollary~\ref{cor-funnycofinal}.
By lemma~\ref{lem-hc}, the composition of the two arrows in the top row is highly connected.
It follows that the map (\ref{eqn-sentry8})
is split injective on homotopy groups or homotopy sets in a certain range. But the second arrow in the top row induces a split injection on
all homotopy groups or homotopy sets, too, as the right-hand part of the diagram shows. \qed

\section{Occupants near the gate} \label{sec-leave}
Imagine a smooth compact manifold
$M$ with corners in the boundary, so that $\partial M=\partial_0M\cup \partial_1M$ as in theorem~\ref{thm-sentriesabsbd}.
(There is no need to assume here that $M$ is a thickening as in definition~\ref{defn-thick}.)
One might ask for a homotopical description of the corner set $\partial\partial_1M=\partial\partial_0M=\partial_0M\cap\partial_1M$ in terms of
the spaces $M\smin S$, where $S$ runs through the finite subsets of $M\smin\partial M$ which are close to
$\partial_0M$. We only make a few steps towards such a description. The
associated definitions turn out to be useful elsewhere \cite{WeissConfHo}.

\subsection{The situation near the gate}
Let $\sP(M\smin\partial_1M)$ be defined as in theorem~\ref{thm-sentriesabsbd}.
Form the twisted arrow poset $\tw(\sP(M\smin\partial_1M))$. (This is a special case
of the twisted arrow category construction which was mentioned in the proof of lemma~\ref{lem-feed1}.)
There is a contravariant functor $\Theta$ from $\tw(\sP(M\smin\partial_1M)$ to spaces given by
\[  \Theta((S,\rho)\le (T,\sigma)):= \big(\textup{closure in $M$ of the collar part of $V(T,\sigma)$}\big)\smin V(S,\rho). \]
This comes with a natural inclusion map
$\partial\partial_1M \to \Theta\big((S,\rho)\le(T,\sigma)\big)$
which in turn induces a map
\begin{equation} \label{eqn-sentry9}
\partial\partial_1M \lra \holim~\Theta\,.
\end{equation}
There is a natural transformation $\Theta \to \Psi\circ F_s$
given by inclusion, where
\[ F_s\co \tw(\sP(M\smin\partial_1M))\lra \sP(M\smin\partial_1M) \]
is the functor \emph{source}. In more detail, for an object $(S,\rho)\le(T,\sigma)$ of the poset $\tw(\sP(M\smin\partial_1M))$ the space
$\Theta((S,\rho)\le(T,\sigma))$ is obviously contained in $\Psi(S,\rho)=M\smin V(S,\rho)$. Using this we get a commutative diagram
\begin{equation} \label{eqn-corner1}
\begin{split}
\xymatrix@C=12pt{
\partial\partial_1M \ar[rr]^-{\textup{inclusion}} \ar[d]^-{(\ref{eqn-sentry9})} && \partial_1M \ar[d]^-{(\ref{eqn-sentry7})} \\
\holim~\Theta \ar[r]  & \holim~\Psi\circ F_s & \ar[l] \holim~\Psi
}
\end{split}
\end{equation}

\begin{lem} \label{lem-corner} The map from $\holim~\Psi$ to $\holim~\Psi\circ F_s$ in~\emph{(\ref{eqn-corner1})}
is a weak equivalence.
\end{lem}

\proof
Recall that $\holim~\Psi$ was defined as $\Tot(X)$ for a certain cosimplicial space $X$. Namely,
$X_r$ is the space of sections of the fiber bundle on $N_r\sP(M\smin\partial M)$ whose fiber over
$((S_0,\rho_0)\ge \cdots\ge (S_r,\rho_r))$ is $M\smin V(S_r,\rho_r)$.
Let $\beta\co \Delta\to \Delta$ be the functor
$[n]\mapsto [2n+1]$. More precisely, $\Delta$ is the category of totally ordered nonempty finite sets
and order-preserving maps, or the equivalent full subcategory with objects $[n]$ for $n\ge 0$, and
$\beta$ is the functor which takes a totally ordered
nonempty finite set $S$ to $S\sqcup S^\op$ (with the total ordering where $a<b$ whenever $a\in S\subset S\sqcup S^\op$
and $b\in S^\op\subset S\sqcup S^\op$). The inclusions
$S\to S\sqcup S^\op$ define a natural transformation $e\co \id\to \beta$.
The cosimplicial space $X\circ\beta$ is Reedy fibrant, by the same argument which we used to show that
$X$ is Reedy fibrant. We need to show that the map $e_*\co \Tot(X)\to \Tot(X\circ\beta)$ is a weak equivalence.
Since both $X$ and $X\circ\beta$ are Reedy fibrant, we can use the easier variant $\Tots$ of $\Tot$
where only the (co)face operators are used; but we continue to view $X$ and $X\circ\beta$ as cosimplicial
spaces. In this setting a more general statement can be made: \emph{if $Y$ is any cosimplicial
space, then $e_*\co \Tots(Y)\to \Tots(Y\circ\beta)$ is a weak equivalence.} To show this,
we use standard resolution procedures
and note that $\Tots$ preserves degreewise weak equivalences.
Therefore we may assume that $Y$ is a homotopy inverse limit of cosimplicial spaces
having the form
\begin{equation}  \label{eqn-cosirep} [r] \mapsto \map(\hom_\Delta([r],[t]),Z) \end{equation}
for some $[t]$ in $\Delta$ and a space $Z$. Since $\Tots$ commutes with such homotopy
inverse limits, it suffices to show that $e_*\co \Tots(Y)\to \Tots(Y\circ\beta)$
is a weak equivalence when $Y$ has the
form~(\ref{eqn-cosirep}). In that case $\Tots(Y)$ is just the space of maps from the geometric
realization of the semi-simplicial set
\begin{equation} \label{eqn-ssc1} [r]\mapsto \hom_\Delta([r],[t]) \end{equation}
to $Z$. Similarly $\Tots(Y\circ\beta)$ is the space of maps from the geometric
realization of the semi-simplicial set
\begin{equation} \label{eqn-ssc2} [r]\mapsto \hom_\Delta([2r+1],[t]) \end{equation}
to $Z$. Now it is enough to show that the geometric realizations of~(\ref{eqn-ssc1})
and~(\ref{eqn-ssc2}) are both contractible. For that we may pretend or observe that
both~(\ref{eqn-ssc1}) and~(\ref{eqn-ssc2}) are actually simplicial sets and realize them as such.
The result is in one case a standard geometric $t$-simplex. In the other case it is an edgewise
subdivided $t$-simplex. \qed

\subsection{A locality statement}
We introduce another topological poset $\sQ=\sQ(M\smin\partial_1M)$
which is closely related to $\tw(\sP)=\tw(\sP(M\smin\partial_1M))$. Intuitively, $\sQ$
is the quotient poset obtained from $\tw(\sP)$ by forcing a morphism
in $\tw(\sP)$ to be an equality (in $\sQ$) if the functor $\Theta$ takes it to an identity map of spaces.
More formally, an element of $\sQ$ is an element
\[  \big((S,\rho) \le (T,\sigma)\big) \]
of $\tw(\sP)$ where $T=\emptyset$, so that $V(T,\sigma)$ is nothing but a collar. For two such pairs
\[  x=\big((S,\rho) \le (T,\sigma)\big), \qquad\quad x'=\big((S',\rho') \le (T',\sigma')\big) \]
where $T=T'=\emptyset$, we say that $x\le x'$ if $\Theta(x')\subset \Theta(x)$.
As a result, $\sQ$ is not a sub-poset of $\tw(\sP)$. Instead
there is a functor (continuous map of posets) $K$ from $\tw(\sP)$ to $\sQ$ which takes
\[  \big((S,\rho) \le (T,\sigma)\big) \]
in $\tw(\sP)$ to $\big((S_0,\rho_0)\le (T_0,\sigma_0)\big)$ where $T_0=\emptyset$ and $\sigma_0$ is the appropriate restriction
of $\sigma$, while $S_0$ is the part of $S$ which is contained in the collar part of $V(T,\sigma)$ and $\rho_0$ is
the appropriate restriction of $\rho$. By construction we have
\[  \Theta=\Theta_1\circ K \]
for a functor $\Theta_1$ from $\sQ$ to spaces. This gives us a prolongation map
\[ \holim~\Theta_1 \lra \holim~\Theta_1\circ K=\holim~\Theta \]
and leads to a commutative diagram
\begin{equation} \label{eqn-corner2}
\begin{split}
\xymatrix@C=17pt{
\partial\partial_1M \ar[rrr]^-{\textup{inclusion}} \ar[d] &&& \partial_1M \ar[d]^-{(\ref{eqn-sentry7})} \\
\holim~\Theta_1 \ar[r] & \holim~\Theta \ar[r]  & \holim~\Psi\circ F_s & \ar[l]_-\simeq \holim~\Psi
}
\end{split}
\end{equation}
which refines diagram~(\ref{eqn-corner1}). Namely, deletion of the vertex $\holim~\Theta_1$ and composition of the
two incident arrows recovers diagram~(\ref{eqn-corner1}).

\medskip
In an older version of this article, it was claimed that the map from $\holim~\Theta_1$ to $\holim~\Theta$ in diagram~(\ref{eqn-corner2})
is always a weak equivalence. The proof had many flaws. (If a more realistic proof of this comes to light, it might go into the next
edition of \cite{WeissConfHo}.) There was also a conjecture which is best restated as saying that
the map from $\partial\partial_1M$ to $\holim~\Theta_1$ in diagram~(\ref{eqn-corner2}) is a weak equivalence under some fairly
restrictive conditions on $\partial_0M$.

\smallskip
In an effort to practice restraint we finish with the following easy but important observation:
the left-hand column in diagram~(\ref{eqn-corner2}) depends only on an arbitrarily
small open neighborhood $U$ of $\partial_0M$ in $M$. To make a precise statement,
let $\partial_1U=U\cap\partial_1M$ and $\partial_0U=\partial_0M$
and define $\sP(U\smin\partial_1U)$ by analogy with $\sP(M\smin\partial_1M)$. Alternatively define it as a full
topological sub-poset of $\sP(M\smin\partial_1M)$, consisting of the objects
$(S,\rho)$ for which the closure of $V(S,\rho)$ in $M$ is contained in $U$. Let $\Theta_{1,U}$ be the restriction
of $\Theta_1$ to $\tw(\sP(U\smin\partial_1U))$.

\begin{prop} The forgetful projection $\holim~\Theta_1\lra \holim~\Theta_{1,U}$
is a weak equivalence.
\end{prop}

\proof The inclusion of the nerve of
$\sQ_U=\sQ(U\smin\partial_1U)$ in the nerve of $\sQ=\sQ(M\smin\partial_1M)$
is a degreewise weak equivalence of simplicial spaces. \qed
\[ *** \]

\medskip
{\bf Remark.} An earlier version of this paper was circulated with M.W. as the sole author.
But S.T. found so many inaccuracies in it that we decided to call the revised version a joint paper
by S.T. and M.W.

\bibliographystyle{amsplain}

\end{document}